\documentclass{amsart}
\usepackage{amsfonts}
\usepackage{amssymb}
\usepackage{xypic}

\newcommand{\cc}{\negmedspace:\negmedspace}
\newcommand{\pomg}{\Omega^{(a)}}
\newcommand{\psig}{\sigma^{(a)}}

\newcommand{\Z}{{\mathbb Z}}
\newcommand{\C}{{\mathbb C}}
\renewcommand{\P}{{\mathbb P}}

\newcommand{\QH}{\text{\sl QH}}

\newcommand{\SL}{\operatorname{SL}}
\newcommand{\Gr}{\operatorname{Gr}}
\newcommand{\Fl}{\operatorname{F\ell}}

\newcommand{\bull}{{\sssize \bullet}}

\newcommand{\Schub}{{\mathfrak{S}}}
\newcommand{\VV}{{\mathcal V}}
\newcommand{\QQ}{{\mathcal Q}}

\newcommand{\gw}[2]{\langle #1 \rangle^{\mbox{}}_{#2}}
\newcommand{\xra}{\xrightarrow}
\newcommand{\til}{\widetilde}
\newcommand{\wbar}{\overline}

\newtheorem{lemma}{Lemma}
\newtheorem{prop}{Proposition}
\newtheorem{thm}{Theorem}

\newtheorem{defn}{Definition}
\theoremstyle{definition}
\newtheorem{example}{Example}

\newcommand{\reflemma}[1]{Lemma~\ref{#1}}
\newcommand{\refthm}[1]{Theorem~\ref{#1}}

\newcommand{\refprop}[1]{Proposition~\ref{#1}}

\newcommand{\refdefn}[1]{Definition~\ref{#1}}

\newcommand{\comment}[1]{}
\newcommand{\noprt}[1]{}

\newenvironment{mathenum}{\begin{enumerate}}{\end{enumerate}}

\newenvironment{romenum}{\begin{enumerate}}{\end{enumerate}}

\begin{document}

\title{Quantum cohomology of partial flag manifolds}
\author{Anders Skovsted Buch}
\date{\today}
\address{Matematisk Institut, Aarhus Universitet, Ny Munkegade, 8000
  {\AA}rhus C, Denmark}
\email{abuch@imf.au.dk}
\maketitle

\section{Introduction}

The (small) quantum cohomology ring of a partial flag variety
$\SL_n(\C)/P$ is a deformation of the usual cohomology ring.  The
structure constants are the three-point, genus zero Gromov-Witten
invariants, which count the number of rational curves meeting three
general Schubert varieties.  The remarkable fact that this ring is
associative \cite{ruan.tian:mathematical,
  kontsevich.manin:gromov-witten} makes it possible to use the
associativity relations to compute Gromov-Witten invariants.

The usual approach for understanding this ring consists of proving a
{\em presentation\/} for the ring \cite{witten:verlinde,
  siebert.tian:on*1, givental.kim:quantum, kim:quantum*2,
  ciocan-fontanine:quantum*1, astashkevich.sadov:quantum, kim:on*3},
together with a {\em quantum Giambelli formula\/} which expresses the
Schubert classes as polynomials in the generators
\cite{bertram:quantum, fomin.gelfand.ea:quantum, ciocan-fontanine:on}.
This information determines the ring as well as all the Gromov-Witten
invariants it encodes.  In addition, a {\em quantum Pieri formula\/}
is known for the multiplication by special Schubert classes
\cite{bertram:quantum, ciocan-fontanine:on, postnikov:on*11}.  These
are the Chern classes of the tautological bundles, and represent the
special Schubert varieties defined by a single Schubert condition.
Since the special Schubert classes generate the quantum ring, the
quantum Pieri formula also determines this ring and its Gromov-Witten
invariants.

The purpose of this paper is to give elementary proofs of the above
structure theorems for the quantum ring of a partial flag variety.  We
do this by proving Ciocan-Fontanine's general quantum Pieri formula
\cite{ciocan-fontanine:on} and by deriving the other results from this
formula.  The quantum Pieri formula is proved by explicitly solving
the underlying Gromov-Witten problem.  That is, given three general
Schubert varieties, one of which is special, we construct the unique
rational curve (of adequate multidegree) meeting these varieties, or
prove that none exist.  We then rely on Ciocan-Fontanine's proof that
the presentation of the quantum ring is a consequence of the quantum
Pieri formula, and give an another argument that the quantum Giambelli
formula is also a consequence.

The original proofs of the quantum formulas relied on intersection
theory on hyperquot schemes.  In the present paper, these techniques
have been replaced with classical Schubert calculus applied to partial
flags called the kernel and span of a curve \cite{buch:quantum,
  buch:direct}.  In particular, if we grant the associativity of
quantum cohomology, we make no use of moduli spaces in this paper.


We also investigate how the theory can best be used to compute
Gromov-Witten invariants.  To this end, we give algorithms for
computing quantum Schubert polynomials and Gromov-Witten invariants.
Despite their simplicity, these algorithms in our experience give an
efficient method for computing in the quantum ring.

In section \ref{sec:cohom} we set up notation and recall the structure
of the usual cohomology ring of a partial flag variety.  We
furthermore give the algorithm for computing quantum Schubert
polynomials (although it is stated for the usual Schubert
polynomials).  In section \ref{sec:quantum} we recall the definition
of the quantum ring, state the quantum Pieri formula, and use it to
derive the remaining results.  We finish this section by explaining
the algorithm for computing Gromov-Witten invariants.  In section
\ref{sec:combin} we prove some combinatorial lemmas relating to the
the classical and quantum Pieri formulas.  Section \ref{sec:tools}
contains geometric tools for handling curves in partial flag
varieties.  These combinatorial and geometric tools are finally used
to prove the quantum Pieri formula in section \ref{sec:proof}.

We thank Sergey Fomin for showing us a very slick proof of
\reflemma{lemma:pieriseq}.  We also thank Ionu{\c{t}} Ciocan-Fontanine
for helpful comments.

\noprt{Spelled wrooooong.}



\section{Cohomology of flag varieties}
\label{sec:cohom}

\subsection{Cohomology}

Set $E = \C^n$.  Given a strictly increasing sequence of integers
$(a_1 < a_2 < \dots < a_k)$ with $a_1 > 0$ and $a_k < n$, we let
$\Fl(a;E)$ be the variety of partial flags $V_{1} \subset V_{2}
\subset \dots \subset V_{k} \subset E$ such that $\dim(V_{i}) =
a_i$ for all $i$.  For convenience we set $a_0 = 0$ and $a_{k+1} =
n$.  The dimension of $\Fl(a;E)$ is equal to $\sum_{i=1}^k
a_i(a_{i+1}-a_i)$.

Let $S_n$ be the group of permutations of $n$ elements.  The Schubert
varieties in $\Fl(a;E)$ are indexed by the set $S_n/W_a$, where $W_a
\subset S_n$ is the subgroup generated by the simple transpositions
$s_i = (i,i+1)$ for $i \not \in \{a_1,\dots,a_k\}$.  Let $S_n(a)
\subset S_n$ denote the set of permutations whose descent positions
are contained in the set $\{a_1,a_2,\dots,a_k\}$.  These permutations
are the shortest representatives for the elements in $S_n/W_a$.  Given
a fixed full flag $F_1 \subset F_2 \subset \dots \subset F_{n-1}
\subset E$ and a permutation $w \in S_n(a)$, define the Schubert
variety
\begin{equation*}
  \Omega^{(a)}_w(F_\bull) = \{ V_\bull \in \Fl(a;E) \mid
  \dim(V_{i} \cap F_p) \geq 
  \#\{t \leq a_i : w(t) > n-p \} ~\forall i,p \} \,.
\end{equation*}
The codimension of this variety is equal to the length $\ell(w)$ of
the permutation $w \in S_n(a)$.

We let $\Omega^{(a)}_w$ denote the fundamental class of
$\Omega^{(a)}_w(F_\bull)$ in the cohomology ring $H^*(\Fl(a;E)) =
H^*(\Fl(a;E); \Z)$.  The Schubert classes $\Omega^{(a)}_w$ form a
basis for this ring, for all $w \in S_n(a)$.  The Schubert class
Poincar{\'e} dual to $\Omega^{(a)}_w$ is the class $\Omega^{(a)}_{w_0 w
  w_a}$ where $w_0 = n\dots 2\, 1$ is the longest permutation in
$S_n$, and $w_a$ is the longest permutation in the subgroup $W_a
\subset S_n$, {i.e.\ }$w_a(j) = a_i + a_{i+1} + 1 - j$ for $a_i < j
\leq a_{i+1}$.

\subsection{Pieri's formula}
\label{sec:pieri}

The Pieri formula gives a rule for multiplying with the Chern classes
of the tautological bundles on $\Fl(a;E)$
\cite{lascoux.schutzenberger:polynomes, sottile:pieris*1}.  Let
$t_{ij}$ denote the transposition interchanging $i$ and $j$.

\begin{defn} \label{defn:arrow}
  Let $1 \leq r \leq m \leq n-1$ be integers and consider the cyclic
  permutation $\alpha = s_r s_{r+1} \cdots s_m \in S_n$ of length
  $\ell = m-r+1$.  For permutations $u$ and $w$ we write $u
  \xra{\alpha} w$ if there exist integers $b_1,\dots,b_\ell$ and
  $c_1,\dots,c_\ell$ such that

\begin{enumerate}
\item $b_i \leq m < c_i$ for all $1 \leq i \leq \ell$;
\item $w = u t_{b_1 c_1} \dots t_{b_\ell c_\ell}$;
\item $\ell(u t_{b_1 c_1} \dots t_{b_i c_i}) = \ell(u) + i$ for all $1
  \leq i \leq \ell$; and
\item the integers $b_1,\dots,b_\ell$ are distinct.
\end{enumerate}
\end{defn}

If $m=a_j$ for some $j$ then $\alpha$ belongs to $S_n(a)$ and
corresponds to the special Schubert variety $\pomg_\alpha(F_\bull)$ of
points $V_\bull \in \Fl(a;E)$ such that $\dim(V_j \cap F_{n-r}) \geq
\ell$.  Its Schubert class is given by $\pomg_\alpha = (-1)^\ell
c_\ell(\VV_j) \in H^*(\Fl(a;E))$, where $\VV_{1} \subset \dots \subset
\VV_{k} \subset E$ denotes the tautological flag on $\Fl(a;E)$.  The
Pieri formula states that for any permutation $u \in S_n(a)$ we have
\begin{equation} \label{eqn:pieri}
  \pomg_\alpha \cdot \pomg_u = \sum_{u \xra{\alpha} w} \pomg_w \,.
\end{equation}


\subsection{Presentation}

We let $\Fl(E) = \Fl(1,2,\dots,n-1; E)$ denote the full flag variety
of $E$, and we denote its Schubert varieties and Schubert classes by
$\Omega_w(F_\bull)$ and $\Omega_w$, respectively, for $w \in S_n$.
The cohomology ring of $\Fl(E)$ has the presentation
\[ H^*(\Fl(E)) = \Z[x_1,\dots,x_n]/(e^n_1,\dots,e^n_n) \]
where $e^n_i = e_i(x_1,\dots,x_n)$ is the $i$th elementary symmetric
polynomial in $n$ variables.  This presentation maps $x_i$ to the
class $\Omega_{s_i} - \Omega_{s_{i-1}}$, which is identical to the
Chern class $-c_1(\VV_i/\VV_{i-1})$.

In this presentation the Schubert class $\Omega_w$ is represented by
the Schubert polynomial $\Schub_w = \Schub_w(x_1,\dots,x_{n-1})$ of
Lascoux and Sch{\"u}tzenberger
\cite{lascoux.schutzenberger:polynomes}.  It is defined as follows.
If $w = w_0$ is the longest permutation in $S_n$, then we set
\[ \Schub_{w_0} = x_1^{n-1} x_2^{n-2} \cdots x_{n-1} \,. \]
Otherwise we can find a simple transposition $s_i \in S_n$ such that
$\ell(w s_i) = \ell(w) + 1$.  In this case we define
\[ \Schub_w =
   \frac{\Schub_{w s_i}(x_1,\dots,x_i,x_{i+1},\dots,x_n) -
   \Schub_{w s_i}(x_1,\dots,x_{i+1},x_i,\dots,x_n)}{x_i - x_{i+1}} \,.
\]
An important property of these polynomials is that they multiply with
the same structure constants as those of the Schubert classes they
represent.  In particular, the Pieri formula (\ref{eqn:pieri}) also
holds as an identity of Schubert polynomials.

The ring $H^*(\Fl(a;E))$ is isomorphic to the subring of
$\Z[x_1,\dots,x_n]/(e^n_1,\dots,e^n_n)$ generated by the elementary
symmetric polynomials $y^p_i = e_i(x_{a_{p-1}+1},\dots,x_{a_p})$ for
$1 \leq p \leq k+1$ and $1 \leq i \leq a_p-a_{p-1}$.  Notice that
$e^n_j = \sum y^1_{i_1} \dots y^{k+1}_{i_{k+1}}$ where the sum is over
all sequences $(i_1, \dots, i_{k+1})$ such that $0 \leq i_p \leq
a_p-a_{p-1}$ and $\sum i_p = j$.  We therefore get the direct
presentation
\begin{equation*}
  H^*(\Fl(a;E)) = \Z[y^1_1,\dots,y^1_{a_1},\, y^2_1,\dots,y^2_{a_2-a_1},\,
  \dots,\, y^{k+1}_1,\dots,y^{k+1}_{n-a_k}] / (e^n_1, \dots, e^n_n)
\end{equation*}
which maps $y^p_i$ to $(-1)^i c_i(\VV_{p}/\VV_{{p-1}}) = c_i(\QQ_p)$, where
$\QQ_p$ is the dual of the bundle $\VV_p/\VV_{p-1}$.

\subsection{An algorithm for Schubert polynomials}
\label{sec:algschub}

If $w \in S_n(a)$ then the Schubert polynomial $\Schub_w$ is symmetric
in each interval of variables $x_{a_{p-1}+1},\dots,x_{a_p}$, so
$\Schub_w$ can be written as a polynomial in the variables $y^p_i$.
This gives a representative of the Schubert class $\Omega^{(a)}_w$ in
the above presentation for $H^*(\Fl(a;E))$.

We will here give a simple method for expressing a Schubert polynomial
$\Schub_w$ for $w \in S_n(a)$ as an integral linear combination
\begin{equation} \label{eqn:fgpexppart}
  \Schub_w = \sum c_{i_{a_1},\dots,i_{n-1}} \, 
  e^{(a)}_{i_{a_1},\dots,i_{n-1}}
\end{equation}
of products of the form
\begin{equation*}
  e^{(a)}_{i_{a_1},\dots,i_{n-1}} = 
  \prod_{p=1}^k \; \prod_{r=a_p}^{a_{p+1}-1} e_{i_r}(x_1,\dots,x_{a_p})
\end{equation*}
for sequences $(i_{a_1},\dots,i_{n-1})$ such that for $a_p \leq r <
a_{p+1}$ we have $0 \leq i_r \leq a_p$.  In fact, if we demand that
$i_{a_p} \leq i_{a_p+1} \leq \dots \leq i_{a_{p+1}-1}$ for all $p$,
then the polynomials $e^{(a)}_{i_{a_1},\dots,i_{n-1}}$ are linearly
independent, so the obtained coefficients $c_{i_{a_1},\dots,i_{n-1}}$
are uniquely determined integers.

Schubert polynomials in the form (\ref{eqn:fgpexppart}) were used by
Fomin, Gelfand, and Postnikov \cite{fomin.gelfand.ea:quantum} and by
Ciocan-Fontanine \cite{ciocan-fontanine:on} to define quantum Schubert
polynomials.  This application will be explained in \S
\ref{sec:qgiambelli}.  Notice that the expression
(\ref{eqn:fgpexppart}) may easily be converted to an expression for
$\Schub_w$ in the $y^p_i$-variables, thus giving the representative of
the class $\Omega^{(a)}_w$ in the presentation for $H^*(\Fl(a;E))$.


The polynomial $\Schub_w$ can be expressed in the form
(\ref{eqn:fgpexppart}) as follows.  Choose $p \leq k$ maximal such
that $w(a_p+1) \neq a_p+1$, and define $u \in S_n(a)$ by
\begin{equation*}
 u(i) = \begin{cases}
   w(i)     & \text{if $i \leq a_p$ and $w(i) < w(a_p+1)$} \\
   w(i)-1   & \text{if $i \leq a_p$ and $w(i) > w(a_p+1)$} \\
   w(i+1)-1 & \text{if $i > a_p$.} \\
 \end{cases}
\end{equation*}
Set $\alpha = s_{w(a_p+1)} \cdots s_{a_p-1} s_{a_p} \in S_n(a)$.  Then
we have $\Schub_\alpha = e^{a_p}_{a_p+1-w(a_p+1)}$.  We claim that the
identity
\begin{equation} \label{eqn:partschub}
  \Schub_w = \Schub_u \cdot e^{a_p}_{a_p+1 - w(a_p+1)}
  - \sum_{u \xra{\alpha} v \neq w} \Schub_v
\end{equation}
can be used recursively to obtain the required expansion of
$\Schub_w$.

Notice that since $w \in S_n(a)$ we automatically have $w(a_p+1) <
a_p+1$.  The identity (\ref{eqn:partschub}) is true by the Pieri
formula because $u \xra{\alpha} w$.  We must show that the recursive
process terminates and that the resulting expression for $\Schub_w$
has the required form (\ref{eqn:fgpexppart}).

For $r < n$ we let $S_r \subset S_n$ denote the subgroup of
permutations fixing the set $\{r+1,\dots,n\}$.  Choose $r$ minimal
such that $w \in S_r$.  Then $r \leq a_{p+1}$ and $u \in S_{r-1}$.
Suppose $\Schub_v$ occurs in the product $\Schub_u \cdot
\Schub_\alpha$.  Then \reflemma{lemma:sketch} of section
\ref{sec:combin}, with $m=a_p$, implies that $w(i) \leq v(i) \leq
u(i)$ for $i \geq a_p+2$.  Now it is immediate from
\refdefn{defn:arrow} that $\sum_{i=a_p+1}^n (u(i)-v(i)) \geq
\ell(\alpha)$.  Since $\ell(\alpha) = \sum_{i=a_p+1}^n (u(i)-w(i))$ we
conclude that $v(a_p+1) \leq w(a_p+1)$, and if equality holds then $v
= w$.  Since we also have $v \in S_r$, this immediately implies
termination.  The resulting expression for $\Schub_w$ is of the form
(\ref{eqn:fgpexppart}) by induction on $r$.

\begin{example}
For $n=7$ and $a = (2,4)$ we get
\[ \Schub_{1536247} = \Schub_{1425367} \cdot e^4_3 - \Schub_{2436157}
   = (e^2_1 e^4_2 - e^4_3) e^4_3 - e^2_1 e^4_1 e^4_4
\]
\end{example}

When using equation (\ref{eqn:partschub}) in real life, it is
essential to remember the Schubert polynomials which have already been
calculated in the recursive process, since otherwise the calculation
of such polynomials can be repeated multiple times.  However, when
this precaution is taken, the algorithm performs well. 
For alternative formulas for Schubert polynomials with relevance to
partial flag varieties we refer to \cite{buch.kresch.ea:schubert}.


\noprt{Spelled wrooooong.}


\section{Quantum cohomology of flag varieties}
\label{sec:quantum}

\subsection{Gromov-Witten invariants}

A rational curve in $\Fl(a;E)$ is the image of a regular map $\P^1 \to
\Fl(a;E)$.  (We will tolerate that a rational curve can be a point
according to this definition.)  The multidegree of a curve $C \subset
\Fl(a;E)$ is the sequence $d = (d_1,\dots,d_k)$ where $d_i$ is the
number of points in the intersection $C \cap \pomg_{s_{a_i}}(F_\bull)$
for any general flag $F_\bull$ of $E$.  Thus, if $C$ is not a point
then the cohomology class of $C$ is equal to $\sum_{i=1}^k d_i \,
\pomg_{w_0 s_{a_i} w_a}$.

Given $u,v,w \in S_n(a)$ and a multidegree $d$ such that $\ell(u) +
\ell(v) + \ell(w) = \dim \Fl(a;E) + \sum_{i=1}^k
d_i(a_{i+1}-a_{i-1})$, the (three-point, genus zero) Gromov-Witten
invariant $\gw{\pomg_u,\pomg_v,\pomg_w}{d}$ is
defined to be the number of rational curves in $\Fl(a;E)$ of
multidegree $d$ meeting all of the Schubert varieties
$\pomg_u(F_\bull)$, $\pomg_v(G_\bull)$, and
$\pomg_w(H_\bull)$ for fixed flags $F_\bull$, $G_\bull$,
$H_\bull$ in general position.  When $\ell(u) + \ell(v) + \ell(w) \neq
\dim \Fl(a;E) + \sum d_i(a_{i+1}-a_{i-1})$ we set
$\gw{\pomg_u,\pomg_v,\pomg_w}{d} = 0$.

Let $q_1, \dots, q_k$ be independent variables and write $\Z[q] =
\Z[q_1,\dots,q_k]$.  The (small) quantum cohomology ring of $\Fl(a;E)$
is a $\Z[q]$-algebra, which as a $\Z[q]$-module is free with a basis
of quantum Schubert classes $\psig_w$:
\begin{equation*}
  \QH^*(\Fl(a;E)) = \bigoplus_{w \in S_n(a)} \Z[q] \, \psig_w \,.
\end{equation*}
Multiplication is defined by the formula
\begin{equation} \label{eqn:qprod}
  \psig_u \cdot \psig_v = \sum_{w,d}
  \gw{\pomg_u,\pomg_v,\pomg_{w_0 w w_a}}{d} \,
  q^d \psig_w
\end{equation}
where the sum is over all $w \in S_n(a)$ and multidegrees $d$, and
$q^d = q_1^{d_1} q_2^{d_2} \cdots q_k^{d_k}$.

It is a non-trivial fact that this product is associative
\cite{ruan.tian:mathematical, kontsevich.manin:gromov-witten,
  fulton.pandharipande:notes}.  The ring $\QH^*(\Fl(a;E))$ has a
natural grading, where the degree of $\psig_w$ is the length
$\ell(w)$, while each variable $q_i$ has degree $a_{i+1}-a_{i-1}$.  If
we set $q_i = 0$ for each $i$, we recover the usual cohomology ring
$H^*(\Fl(a;E))$.

\subsection{The quantum Pieri formula}
\label{sec:qpieri}

The central result about the structure of the quantum ring
$\QH^*(\Fl(a;E))$ is the quantum Pieri formula of Ciocan-Fontanine
\cite{ciocan-fontanine:on}.  This result generalizes the quantum Pieri
formula for Grassmannians \cite{bertram:quantum} and the quantum
Monk's formula for full flag varieties
\cite{fomin.gelfand.ea:quantum}.  In the case of full flag varieties,
Postnikov has given an equivalent but simpler statement of the quantum
Pieri formula, as well as a combinatorial proof based on the quantum
Monk's formula \cite{postnikov:on*11}.  We will give an elementary
geometric proof of Ciocan-Fontanine's result in the last section.

We will call a sequence $d = (d_1,\dots,d_k)$ of non-negative integers
for a {\em Pieri sequence with maximum at position $j$}, if
$(d_1,\dots,d_j)$ is weakly increasing, $(d_j,\dots,d_k)$ is weakly
decreasing, and if we set $d_0 = d_{k+1} = 0$ then $|d_i - d_{i+1}|
\leq 1$ for $0 \leq i \leq k$.  Given such a sequence, we set
$\gamma_d = \tau_1 \tau_2 \cdots \tau_k \in S_n$ where $\tau_i$ is the
permutation which interchanges the intervals $[a_i-d_i+1,a_i]$ and
$[a_i+1,a_{i+1}]$.  In other words, $\tau_i$ is defined by
\[ \tau_i(p) = \begin{cases}
  p+a_{i+1}-a_i & \text{if $a_i-d_i < p \leq a_i$}, \\
  p-d_i & \text{if $a_i < p \leq a_{i+1}$}, \\
  p & \text{otherwise.}
\end{cases} \]

\begin{thm}[Quantum Pieri formula \cite{ciocan-fontanine:on}] 
\label{thm:partqpieri}
Let $\alpha = s_r s_{r+1} \cdots s_{a_j}$ and $u \in
S_n(a)$ be permutations.  Then
\begin{equation*}
\psig_\alpha \cdot \psig_u = \sum q^d \psig_w
\end{equation*}
where the sum is over all Pieri sequences $d$ with maximum at position
$j$ and permutations $w \in S_n(a)$ such that (i) $\ell(u \gamma_d) =
\ell(u) - \ell(\gamma_d)$; (ii) $\ell(w w_a \gamma_d) = \ell(w w_a) +
\ell(\gamma_d)$; and (iii) $u \gamma_d \xra{\wbar \alpha} w w_a
\gamma_d w_b$ where $b = a - d = (a_1-d_1,\dots,a_k-d_k)$ and $\wbar
\alpha = s_r s_{r+1} \cdots s_{b_j}$.
\end{thm}

An equivalent symmetric version of this theorem is given in section
\ref{sec:proof}.  Notice that condition (iii) implicitly implies that
$d_j \leq \ell(\alpha)$.  \comment{FIXME: Are (i) and (ii) equivalent
  with $\ell(w) = \ell(\alpha) + \ell(u) - \sum
  d_i(a_{i+1}-a_{i-1})$?}

Given a Pieri sequence $d$ with maximum at position $j$, set $h_p =
\min \{i : d_i = p\}$ and $l_p = \max \{ i : d_i = p \}$ for each $1
\leq p \leq d_j$.  With this notation we have $\gamma_d(b_{l_p+1}) =
a_{h_p}$, while $\gamma_d(i) = i+p$ if $b_{h_p} < i \leq b_{h_{p+1}}$
or if $b_{l_{p+1}+1} < i < b_{l_p+1}$.  It follows that $\ell(u
\gamma_d) = \ell(u) - \ell(\gamma_d)$ if and only if $u(a_{h_p}) >
u(i)$ for all $p$ and $a_{h_p} < i \leq a_{l_p+1}$ (cf. \cite[Remark
3.2 (ii)]{ciocan-fontanine:on}.)

\begin{example}
  Let $\alpha = s_2 s_3 s_4$ and $u = 3715246$.  We will compute the
  product $\psig_\alpha \cdot \psig_u$ in the ring $\QH^*
  \Fl(2,4;\C^7)$.  First observe that the Pieri sequences $d$ with
  maximum at position $2$ such that $\ell(u \gamma_d) = \ell(u) -
  \ell(\gamma_d)$ are $(0,0)$ and $(1,1)$.  The first of these
  contributes with $\sum_{u \xra{\alpha} w} \psig_w =
  \psig_{4726135}$.  For $d=(1,1)$ we have $\wbar \alpha = s_2 s_3$,
  and $u \gamma_d = 3152467 \xra{\wbar \alpha} v$ when $v$ is one of
  the permutations $4251367$, $3261457$, and $4162357$.  The first two
  of these satisfy $\ell(v w_b \gamma_d^{-1}) = \ell(v w_b) -
  \ell(\gamma_d)$, and they contribute $q_1 q_2\, \psig_{1425367} +
  q_1 q_2\, \psig_{1326457}$.  In conclusion we have
\[ \psig_\alpha \cdot \psig_u = \psig_{4726135} + 
   q_1 q_2\, \psig_{1425367} + q_1 q_2\, \psig_{1326457} \,.
\]
\end{example}

\subsection{Structure of the quantum ring}
\label{sec:structure}

The presentation of $\QH^* \Fl(a;E)$ is due to Astashkevich and Sado
\cite{astashkevich.sadov:quantum} and Kim \cite{kim:quantum*3,
  kim:on*3} (see also \cite{witten:verlinde, siebert.tian:on*1} for
the Grassmannian case, and \cite{givental.kim:quantum, kim:quantum*2,
  ciocan-fontanine:quantum*1} for the case of full flag varieties.)
In this section we sketch how to recover this presentation from the
quantum Pieri formula.  We follow Ciocan-Fontanine's paper
\cite{ciocan-fontanine:on}.

Let $\phi : H^*(\Fl(a;E)) \to \QH^*(\Fl(a;E))$ be the linear map which
sends each Schubert class $\pomg_w$ to the corresponding quantum
Schubert class $\psig_w$.  The presentation of $\QH^*(\Fl(a;E))$ uses
variables $y^j_i$ and $q_j$, and maps $y^j_i$ to $(-1)^i
\phi(c_i(\VV_j/\VV_{j-1})) = \phi(c_i(\QQ_j))$.

Set $\alpha_{i,j} = s_{a_j-i+1} s_{a_j-i+2} \cdots s_{a_j}$ and
$\beta_{i,j} = s_{a_j+i-1} s_{a_j+i-2} \cdots s_{a_j} \in S_n(a)$.
Then $\Schub_{\alpha_{i,j}} = e_i(x_1,\dots,x_{a_j})$ is the
elementary symmetric polynomial and $\Schub_{\beta_{i,j}} =
h_i(x_1,\dots,x_{a_j})$ is the complete symmetric polynomial in $a_j$
variables.  Using the Pieri formula (\ref{eqn:pieri}), it follows that
for $i \leq a_j - a_{j-1}$ we have (cf.\ \cite[Lemma
3.5]{ciocan-fontanine:on})
\begin{equation*}
  c_i(\QQ_j) 
  = \sum_{p=0}^i (-1)^p\, \pomg_{\alpha_{i-p,j}} \cdot
    \pomg_{\beta_{p,j-1}}
  = \sum_{p=0}^i (-1)^p\, \pomg_{\beta_{p,j-1} \alpha_{i-p,j}} \,.
\end{equation*}
We therefore get
\begin{equation*}
  \phi(c_i(\QQ_j)) = \sum_{p=0}^i 
  (-1)^p\, \psig_{\beta_{p,j-1} \alpha_{i-p,j}} \,.
\end{equation*}

Define {\em quantum elementary symmetric polynomials\/} $E^j_i$ as
follows.  If $i = j = 0$ then set $E^0_0 = 1$.  If $j<0$ or $i<0$ or
$i>a_j$ then set $E^j_i = 0$.  Otherwise, if $0 \leq i \leq a_j \neq
0$ then define inductively
\begin{equation*}
  E^j_i = E^{j-1}_i + \sum_{r=1}^{a_j-a_{j-1}} y^j_r \, E^{j-1}_{i-r} 
  \, - \, (-1)^{a_j-a_{j-1}}\, q_{j-1} E^{j-2}_{i-a_j+a_{j-2}} \,.
\end{equation*}
For example, if $n=7$ and $a = (2,4)$ then we get 
\[\begin{split} E^3_5 &= E^2_5 + y^3_1 E^2_4 + y^3_2 E^2_3 + y^3_3
  E^2_2 - (-1)^3 q_2 E^1_0 \\
&= 0 + y^3_1(y^2_2 y^1_2 - q_1) + y^3_2(y^2_1 y^1_2 + y^2_2 y^1_1) +
y^3_3(y^1_2 + y^2_1 y^1_1 + y^2_2) + q_2 \,.
\end{split}\]

We claim that the quantum ring has the presentation
\[ \QH^*(\Fl(a;E)) = \Z[y,q]/(E^{k+1}_1, \dots, E^{k+1}_n) \]
where each variable $y_i^j$ is mapped to $\phi(c_i(\QQ_j))$ for $1
\leq j \leq k+1$ and $1 \leq i \leq a_j-a_{j-1}$.

More generally, if we replace each $y^j_i$ with $\phi(c_i(\QQ_j))$,
then $E_i^j$ maps to $\psig_{\alpha_{i,j}}$ for $j \leq k$ while
$E_i^{k+1}$ becomes zero.  In fact, since a symmetric functions
calculation shows that this is true in cohomology after setting
$q_j=0$ for all $j$, we only need to determine the $q$-terms which
arise when the sum $\sum_{r=1}^{a_j-a_{j-1}} \phi(c_r(\QQ_j))\cdot
\psig_{\alpha_{i-r,j-1}}$ is expanded.  Here one observes that, if $d$
is a non-zero Pieri sequence with maximum at position $j-1$ and if
$\ell(\beta_{p,j-1} \alpha_{r-p,j} \gamma_d) = \ell(\beta_{p,j-1}
\alpha_{r-p,j}) - \ell(\gamma_d)$, then $r=p=a_j-a_{j-1}$ and $d =
(0,\dots,0,1,0,\dots,0)$ has a single one at position $j-1$.
Furthermore, when $r = a_j-a_{j-1}$ the product $\psig_{\beta_{r,j-1}}
\cdot \psig_{\alpha_{i-r,j-1}}$ contains no $q$-terms for $i <
a_j-a_{j-2}$ while it has exactly one when $i \geq a_j-a_{j-2}$,
namely $q_{j-1} \psig_{\alpha_{i-a_j+a_{j-2}, j-2}}$.  For more
details we refer to \cite[Lemma 3.6]{ciocan-fontanine:on}.


\subsection{The quantum Giambelli formula}
\label{sec:qgiambelli}

For a sequence $(i_{a_1},\dots,i_{n-1})$ such that $0 \leq i_r \leq a_j$
for each $a_j \leq r < a_{j+1}$, set
\begin{equation*}
  E^{(a)}_{i_{a_1},\dots,i_{n-1}} = 
  \prod_{j=1}^k \prod_{r=a_j}^{a_{j+1}-1} E^j_{i_r} \,.
\end{equation*}
Ciocan-Fontanine has given a geometric proof that no $q$-terms occur
in the expansion of the corresponding product of quantum Schubert
classes \cite[Thm.\ 3.14]{ciocan-fontanine:on}.  In other words,
$\phi$ maps the cohomology class represented by
$e^{(a)}_{i_{a_1},\dots,i_{n-1}}$ to the quantum class given by
$E^{(a)}_{i_{a_1},\dots,i_{n-1}}$.  We will here deduce this fact from
\refthm{thm:partqpieri}.

Suppose $u \in S_n(a) \cap S_r$ for some $r < a_{j+1}$.  Then $\psig_u
\cdot \psig_{\alpha_{i,j}} = \phi(\pomg_u \cdot
\pomg_{\alpha_{i,j}})$, which follows because $\ell(u \gamma_d) >
\ell(u) - \ell(\gamma_d)$ for all non-zero Pieri sequences $d$ with
maximum at position $j$.  \reflemma{lemma:sketch} of section
\ref{sec:combin} furthermore implies that all terms $\pomg_w$ in the
product $\pomg_u \cdot \pomg_{\alpha_{i,j}}$ satisfy $w \in S_n(a)
\cap S_{r+1}$.  For each $a_1 \leq r \leq n-1$ we set $\alpha(r) =
\alpha_{i_r,j}$ where $j$ is maximal such that $a_j \leq r$.  By
induction on $r$, the above comments imply that $\psig_{\alpha(a_1)}
\cdot \psig_{\alpha(a_1+1)} \cdots \psig_{\alpha(r)} =
\phi(\pomg_{\alpha(a_1)} \cdot \pomg_{\alpha(a_1+1)} \cdots
\pomg_{\alpha(r)})$ and that all terms $\psig_w$ in this product
satisfy $w \in S_{r+1}$.  In particular, the class given by
$E^{(a)}_{i_{a_1},\dots,i_{n-1}}$ contains no $q$-terms.

Define the {\em partial quantum Schubert polynomial\/} for a
permutation $w \in S_n(a)$ by
\begin{equation*}
  \Schub^q_w = \sum c_{i_{a_1},\dots,i_{n-1}} \, 
  E^{(a)}_{i_{a_1},\dots,i_{n-1}} \,.
\end{equation*}
The coefficients $c_{i_{a_1},\dots,i_{n-1}}$ are defined by
(\ref{eqn:fgpexppart}).  In particular, this definition depends on the
sequence $a$.

By applying $\phi$ to the classes represented by either side of
(\ref{eqn:fgpexppart}), it follows that $\Schub^q_w$ is a
representative for the quantum Schubert class $\psig_w$ in the
presentation of $\QH^*(\Fl(a;E))$.  In other words, $\Schub^q_w$ is a
quantum Giambelli formula.  This result is due to Bertram
\cite{bertram:quantum} for Grassmannians, to Fomin, Gelfand, and
Postnikov for full flag varieties \cite{fomin.gelfand.ea:quantum}, and
to Ciocan-Fontanine in general \cite{ciocan-fontanine:on}.

Notice that the identity (\ref{eqn:partschub}) gives the direct
recursive formula
\begin{equation} \label{eqn:qschubpart}
  \Schub^q_w = \Schub^q_u \cdot E^p_{a_p+1-w(a_p+1)} - 
  \sum_{u \xra{\alpha} v \neq w} \Schub^q_v
\end{equation}
where $p$, $u$, and $\alpha$ are chosen as in \S \ref{sec:algschub}.
For full flag varieties one can alternatively use a quantum version of
the transition formula, which is based on the quantum Monk's formula
(see \cite[\S 8]{fomin.gelfand.ea:quantum} and
\cite[(4.16)]{macdonald:notes}).

\subsection{Computing Gromov-Witten invariants}

By definition of the quantum product (\ref{eqn:qprod}), a
Gromov-Witten invariant $\gw{\pomg_u, \pomg_v, \pomg_w}{d}$ on
$\Fl(a;E)$ can be computed by extracting the coefficient of
$q^d \psig_{w_0 w w_a}$ in the expansion of the product $\psig_u \cdot
\psig_v \in \QH^*(\Fl(a;E))$.

The following method may be used to compute this product.  Start by
expressing the quantum class $\psig_u$ as a polynomial in the classes
$\psig_{\alpha_{ij}}$.  This can be done using equation
(\ref{eqn:qschubpart}).  Then let this polynomial act on the class
$\psig_v$ using the quantum Pieri formula (Thm.\ 
\ref{thm:partqpieri}).  The result is the desired expansion.

Practical experiments indicate that this method is quite efficient.
For example, it vastly outperforms the Gr{\"o}bner basis methods
suggested in \cite{fomin.gelfand.ea:quantum}.  Notice that the roles
of $u$, $v$, and $w$ can be permuted.  Often (but not always) the best
choice is to let the quantum Schubert polynomial for the shortest
permutation act on one of the other quantum Schubert classes.  Notice
also that this method for computing Gromov-Witten invariants does not
make any use of the presentation of the quantum ring.


\begin{example}
  We will compute the number of rational curves in $\Fl(2,4;\C^7)$ of
  multidegree $(2,3)$, which pass through two general points and meet
  the Schubert variety $\pomg_{3417256}(F_\bull)$.  In other words, we
  compute the Gromov-Witten invariant $\gw{\pomg_{6745123},
    \pomg_{6745123}, \pomg_{3417256}}{(2,3)}$.  Using equation
  (\ref{eqn:qschubpart}) we obtain the point class as $\psig_{6745123}
  = (\psig_{\alpha_{2,1}})^2 (\psig_{\alpha_{4,2}})^3$.  Now using the
  quantum Pieri formula repeatedly, we obtain
\[\begin{split} 
  \psig_{6745123} \cdot \psig_{3417256} 
  &= (\psig_{\alpha_{2,1}})^2 (\psig_{\alpha_{4,2}})^3 \cdot
     \psig_{3417256} \\
  &= (\psig_{\alpha_{2,1}})^2 (\psig_{\alpha_{4,2}})^2 \cdot 
     q_2\, \psig_{4512367} \\
  &= (\psig_{\alpha_{2,1}})^2 \psig_{\alpha_{4,2}} \cdot
     q_2\, \psig_{5623147} \\
  &= (\psig_{\alpha_{2,1}})^2 \cdot q_2\, \psig_{6734125} \\
  &= \psig_{\alpha_{2,1}} \cdot 
     (q_1 q_2\, \psig_{3746125} + q_1 q_2^2\, \psig_{1734256}) \\
  &= q_1^2 q_2\, \psig_{3467125} + q_1^2 q_2^2\, \psig_{1436257} +
     q_1^2 q_2^2\, \psig_{1347256} + q_1^2 q_2^3\, \psig_{1234567} \,.
\end{split}\]
The Gromov-Witten invariant of interest is the coefficient to $q_1^2
q_2^3\, \psig_{1234567}$ in this product, so it is equal to one.
\end{example}

\noprt{Spelled wrooooong.}


\section{Combinatorics of the Pieri rule}
\label{sec:combin}

In this section we prove some lemmas concerning the Pieri and quantum
Pieri formulas.  As in \S \ref{sec:pieri} we set $\alpha = s_r s_{r+1}
\cdots s_m \in S_n$ and $\ell = \ell(\alpha)$.

\begin{lemma} \label{lemma:permute}
  Let ${u} \xra{\alpha} {w}$ and let $b$ and $c$ be sequences
  satisfying \refdefn{defn:arrow}.  Suppose $c_i \neq c_{i+1}$.  Then
  we can interchange the $i$'th and $(i+1)$'th indices in $b$ and $c$,
  i.e.\ the sequences $b' = (b_1,\dots,b_{i+1},b_i,\dots,b_\ell)$ and
  $c' = (c_1,\dots,c_{i+1},c_i,\dots,c_\ell)$ also satisfy
  \refdefn{defn:arrow}.
\end{lemma}
\begin{proof}
  The sequences $b'$ and $c'$ clearly satisfy properties (1) and (4).
  Conditions (2) and (3) hold because the transposition $t_{b_i c_i}$
  commutes with $t_{b_{i+1} c_{i+1}}$.
\end{proof}

The following fact has already be used in \S \ref{sec:algschub} and \S
\ref{sec:structure}.

\begin{lemma} \label{lemma:sketch}
  Let ${u} \xra{\alpha} {w}$ and suppose ${u}$ has no descents after
  position $m$.  Then for all $j \geq m+2$ we have ${u}(j-1) < {w}(j)
  \leq {u}(j)$.
\end{lemma}
\begin{proof}
  Let $b$ and $c$ be sequences satisfying \refdefn{defn:arrow}.  By
  \reflemma{lemma:permute} we may assume that $c_1=\dots=c_p=j$ and
  $c_i \neq j$ for $i > p$.  It then follows from property (3) and
  induction on $i$ that each permutation ${u} t_{b_1 c_1} \cdots
  t_{b_i c_i}$ maps $j$ to a value greater than ${u}(j-1)$.
\end{proof}

\begin{lemma} \label{lemma:misc1}
  Let ${u} \xra{\alpha} {w}$ and suppose that for some $l > m$ we have
  ${u}(i) < {u}(l)$ for all $m < i < l$.  Then for all $j \leq m$ such
  that ${u}(j) < {u}(l)$ we have ${w}(j) \leq {u}(l)$.
\end{lemma}
\begin{proof}
  Let $b$ and $c$ be sequences satisfying \refdefn{defn:arrow}.  We
  may assume that ${w}(j) \neq {u}(j)$, so $j = b_p$ for some $p$.  By
  \reflemma{lemma:permute} we may furthermore assume that $c_1 = c_2 =
  \dots = c_p$.
  
  Set $u' = {u} t_{b_1 c_1} \cdots t_{b_{p-1} c_{p-1}}$.  If $c_p \leq
  l$ then ${w}(j) = u'(c_p) \leq {u}(c_p) \leq {u}(l)$.  On the other
  hand, if $c_p > l$ then since $u'(j) < u'(l)$ and $\ell(u' t_{b_p
    c_p}) = \ell(u') + 1$ we must have $u'(c_p) < u'(l)$, so once
  again we get ${w}(j) = u'(c_p) < u'(l) = {u}(l)$, as required.
\end{proof}

If $x_1, \dots, x_p$ are elements of a vector space $E$, we let
$\left< x_1,\dots, x_p \right> \subset E$ denote the linear span of
these vectors.

\begin{lemma} \label{lemma:basis}
  Let $\{ e_1, \dots, e_n \}$ be a basis for a vector space $E$ and
  let ${u}, {w} \in S_n$ be permutations such that ${u}
  \xra{\alpha} {w}$.  Suppose $x_1, \dots, x_n \in E$ are elements
  satisfying the following conditions:
\begin{romenum}
\item If $i \leq m$ and ${u}(i) = {w}(i)$ then $x_i = e_{{u}(i)}$
\item If $i \leq m$ and ${u}(i) \neq {w}(i)$ then $x_i = \lambda_i
  e_{{u}(i)} + \mu_i e_{{w}(i)}$ where $\lambda_i, \mu_i \neq 0$
\item If $i > m$ then $x_i = e_{{u}(i)}$ or $x_i = e_{{w}(i)}$.
\end{romenum}
Then $\{x_1, \dots, x_n\}$ is also a basis for $E$.  The flag $V_\bull
\in \Fl(E)$ given by $V_i = \left< x_1, \dots, x_i \right>$ belongs to
the Schubert variety $\Omega_{u}(F_\bull)$ where $F_\bull$ is defined
by $F_i = \left< e_{n+1-i}, \dots, e_n \right>$.  Furthermore, this
flag $V_\bull$ does not depend on the choices made in (iii).
\end{lemma}
\begin{proof}
  Suppose at first that $x_i = e_{{u}(i)}$ for all $i > m$.  In this
  case we have $F_{n+1-{u}(i)} = F_{n-{u}(i)} \oplus \C x_i$ for all
  $i$, so $\{ x_1, \dots, x_n\}$ is a basis and $V_\bull \in
  \Omega_{u}(F_\bull)$.  It suffices to show that $V_{i-1} \oplus \C
  e_{{u}(i)} = V_{i-1} \oplus \C e_{{w}(i)}$ for each $i > m$.  In
  case $u(i) \neq w(i)$ we let $b$ and $c$ be sequences satisfying
  \refdefn{defn:arrow}, such that for some $p$ we have $c_1 = \dots =
  c_p = i$ and $c_j \neq i$ for $j > p$.  Then the values
  $w(b_1),w(b_2),\dots,w(b_p),w(i)$ agree with
  $u(i),u(b_1),\dots,u(b_{p-1}),u(b_p)$, in the indicated order.  This
  implies that $\left< x_{b_1},\dots,x_{b_p} \right>$ is a subspace of
  $\left< e_{u(b_1)}, \dots, e_{u(b_p)}, e_{u(i)} \right>$.
  Furthermore, (ii) implies that neither $e_{u(i)}$ nor $e_{w(i)}$ is
  contained in this subspace, so $\left< x_{b_1},\dots,x_{b_p},
    e_{u(i)} \right> = \left< x_{b_1},\dots,x_{b_p}, e_{w(i)} \right>
  = \left< e_{u(b_1)}, \dots, e_{u(b_p)}, e_{u(i)} \right>$.  The
  required identity of subspaces follows from this.
\end{proof}

We also need the following characterization of Pieri sequences, which
is equivalent to parts (i) and (ii) of \cite[Lemma
5.2]{ciocan-fontanine:on}.

\begin{lemma} \label{lemma:pieriseq}
  A sequence of non-negative integers $d = (d_1, \dots, d_k)$ is a
  Pieri sequence with maximum at position $j$ if and only if the
  inequality
\[ d_j + \sum_{i=1}^{k-1} d_i d_{i+1} - \sum_{i=1}^k d_i^2 ~\geq~ 0 \]
is satisfied.  In this case the inequality is satisfied with equality.
\end{lemma}
\begin{proof}(Fomin)
  The inequality can be rewritten as $2 d_j \geq \sum_{i=0}^k (d_i -
  d_{i+1})^2$ and the right-hand side of this is estimated from below
  by $\sum |d_i - d_{i+1}| \geq 2 d_j$.
\end{proof}

\noprt{Spelled wrooooong.}


\section{Geometric tools}
\label{sec:tools}

In this section we will give some tools for handling curves in flag
varieties.  It is to convenient to extend the notation for partial
flag varieties to allow weakly increasing sequences of dimensions.  If
$b = (b_1 \leq b_2 \leq \dots \leq b_k)$ is a weakly increasing
sequence with $b_1 \geq 0$ and $b_k \leq n$ we let $\Fl(b;E)$ be the
variety of partial flags $K_{1} \subset K_{2} \subset \dots \subset
K_{k} \subset E$ such that $\dim K_{i} = b_i$ for all $i$.  The
Schubert varieties in $\Fl(b;E)$ are indexed by the set $S_n(b)$ of
permutations whose descent positions are contained in
$\{b_1,\dots,b_k\}$.

Let $b$ be a weakly increasing sequence such that $b_i \leq a_i$ for
each $i$.  Given a Schubert variety $\pomg_w(F_\bull) \subset
\Fl(a;E)$ we will need a description of the set of points $K_\bull \in
\Fl(b;E)$ such that for some $V_\bull \in \pomg_w(F_\bull)$ we have
$K_{i} \subset V_{i}$ for all $i$.

We construct a permutation $\overline w \in S_n(b)$ from $w$ as
follows.  Set $w^{(0)} = w$.  Then for each $1 \leq i \leq k$ we let
$w^{(i)}$ be the permutation obtained from $w^{(i-1)}$ by rearranging
the elements $w^{(i-1)}(b_i+1), \dots, w^{(i-1)}(a_{i+1})$ in
increasing order.  Finally we set $\overline w = w^{(k)}$.  For
example, if $n = 6$, $a = (2,5)$, $b = (1,2)$, and $w =
2\,6\,3\,4\,5\,1$ then $w^{(1)} = 2\,3\,4\,5\,6\,1$ and $\overline w =
2\,3\,1\,4\,5\,6$.  The following result is proved in
\cite{buch:direct}.

\begin{lemma} \label{lemma:kernel}
  The set $\{ K_\bull \in \Fl(b;E) \mid \exists~ V_\bull \in
  \Omega^{(a)}_w(F_\bull) : K_{i} \subset V_{i}
  ~\forall i \}$ is equal to the Schubert variety $\Omega_{\overline
  w}^{(b)}(F_\bull)$ in $\Fl(b;E)$.
\end{lemma}

Our notation is related to Pieri sequences as follows.

\begin{lemma} \label{lemma:gammabar}
  Let $d$ be a Pieri sequence and set $b = a - d =
  (a_1-d_1,\dots,a_k-d_k)$.  Let $u \in S_n(a)$.  Then $\ell(u
  \gamma_d) = \ell(u) - \ell(\gamma_d)$ if and only if $\ell(\wbar u) =
  \ell(u) - \ell(\gamma_d)$ if and only if $u \gamma_d = \wbar u$.  In
  this case we have $u \gamma_d \in S_n(b)$.
\end{lemma}
\begin{proof}
  With the notation of section \ref{sec:qpieri} we have
  $\ell(\gamma_d) = \sum \ell(\tau_i)$.  The lemma follows because
  $\ell(u^{(i-1)} \tau_i) \geq \ell(u^{(i-1)}) - \ell(\tau_i)$ and
  $\ell(u^{(i)}) \geq \ell(u^{(i-1)}) - \ell(\tau_i)$ for all $i$,
  with equality if and only if $u^{(i)} = u^{(i-1)} \tau_i$.
\end{proof}

Now let $C \subset \Fl(a;E)$ be a rational curve of multidegree $d =
(d_1,\dots,d_k)$.  For each $i$ we let $C_i = \rho_i(C) \subset
\Gr(a_i,E)$ be the image of $C$ in the Grassmannian $\Gr(a_i,E)$ by
the projection $\rho_i : \Fl(a;E) \to \Gr(a_i,E)$.  This curve $C_i$
then has a {\em kernel\/} and a {\em span\/} \cite{buch:quantum}.  The
kernel is the largest subspace of $E$ contained in all the
$a_i$-dimensional subspaces of $E$ corresponding to points of $C_i$.
We let $b_i$ be the dimension of this kernel and denote the kernel
itself by $K_{i}$.  It follows from \cite[Lemma 1]{buch:quantum} that
$b_i \geq a_i - d_i$ for each $i$.  The span of $C_i$ is the smallest
subspace of $E$ containing all subspaces given by points of $C_i$.
This span has dimension at most $a_i - d_i$.

The kernels $K_{i}$ form a partial flag $K_\bull \in \Fl(b;E)$ called
the kernel of $C$.  Notice that $K_{i} \subset V_{i}$ for all points
$V_\bull \in C$.  \reflemma{lemma:kernel} therefore implies the
following (cf.\ \cite[Prop.\ 1]{buch:direct}).

\begin{prop} \label{prop:direct}
Let $C \in \Fl(a;E)$ be a rational curve with kernel $K_\bull \in
\Fl(b;E)$.  If $C \cap \pomg_w(F_\bull) \neq \emptyset$ then $K_\bull
\in \Omega^{(b)}_{\overline w}(F_\bull)$.
\end{prop}

\begin{lemma} \label{lemma:minkernel}
  Let $f : \P^1 \to \Gr(m,E)$ be a curve of degree $d$ such that the
  kernel $K$ of $f(\P^1)$ has dimension $m-d$.  Then there are
  elements $x_1,\dots,x_d, y_1, \dots, y_d \in E$ such that $f(s\cc t)
  = K \oplus \left< s x_1 + t y_1, \dots, s x_d + t y_d \right>$ for
  all $(s\cc t) \in \P^1$.
\end{lemma}
\begin{proof}
  Any regular map $f : \P^1 \to \Gr(m,E)$ can be written in the form
  $f(s\cc t) =\linebreak \left< f_1(s\cc t), \dots, f_m(s\cc t) \right>$ for
  regular maps $f_i : \P^1 \to \P(E)$, and furthermore we have $\sum
  \deg(f_i) = \deg(f) = d$.  (To see this, one uses that the pullback
  of the tautological subbundle on $\Gr(m;E)$ splits as a sum of line
  bundles on $\P^1$.)  At least $m-d$ of these maps must have degree
  zero, so we can assume that $f_{d+1}, \dots, f_m$ are constant.
  Since $\left< f_{d+1}, \dots, f_m \right>$ is contained in $K$ and
  these spaces have the same dimension, we conclude that $K = \left<
    f_{d+1}, \dots, f_m \right>$.  This implies that none of the
  functions $f_1, \dots, f_d$ are constant, so they must all have
  degree one.  Thus we can write $f_i(s\cc t) = s x_i + t y_i$ for
  some $x_i, y_i \in E$ for $1 \leq i \leq d$.
\end{proof}


Given a morphism $f : \P^1 \to \Fl(a;E)$ we let $f_i : \P^1 \to
\Gr(a_i,E)$ denote the composition of $f$ with the $i$th projection
$\rho_i : \Fl(a;E) \to \Gr(a_i,E)$.

\begin{lemma} \label{lemma:classify}
  Let $a = (a_1 < a_2 < a_3)$ be a sequence of integers, $0 < a_i <
  n$, and let $f = (f_1,f_2,f_3) : \P^1 \to \Fl(a;E)$ be a regular map
  of multidegree $(d, d+1, d)$ for some integer $d \geq 0$.  Suppose
  that the kernel $K_\bull \in \Fl(b;E)$ of $f(\P^1)$ has dimensions
  given by $b = (a_1-d, a_2-d-1, a_3-d)$.  Suppose also there are
  linearly independent elements $x_1, \dots, x_{d+1}, y_1, \dots,
  y_{d+1} \in E$ such that
\begin{mathenum}
\item $f_1(s\cc t) = K_{1} \oplus \left< s\, x_1 + t\, y_1, \dots, s\,
    x_d + t\, y_d \right>$ for all $(s\cc t) \in \P^1$
\item $f_2(1\cc 0) = K_{2} \oplus \left< x_1, \dots, x_{d+1} \right>$
\item $f_2(0\cc 1) = K_{2} \oplus \left< y_1, \dots, y_{d+1} \right>$
\item $K_{2} \cap \left< x_1,\dots,x_{d+1}, y_1,\dots,y_{d+1} \right> = 0$
\item $K_{3} \cap \left< x_1,\dots,x_d, y_1,\dots,y_d \right> = 0$
\item $x_{d+1}, y_{d+1} \in K_{3}$.
\end{mathenum}
Then there exists a unique $\lambda \in \C^*$ such that
\[ f_2(s\cc t) = K_{2} \oplus \left< s\, x_1 + t\, y_1, \dots, s\, x_d +
  t\, y_d, s\, x_{d+1} + t\, \lambda\, y_{d+1} \right> \,.
\]
\end{lemma}
\begin{proof}
  By \reflemma{lemma:minkernel} we can find elements $x'_1, \dots,
  x'_{d+1}, y'_1, \dots, y'_{d+1} \in E$ such that $f_2(s\cc t) = K_{2}
  \oplus \left< s\, x'_1 + t\, y'_1, \dots, s\, x'_{d+1} + t\,
    y'_{d+1} \right>$.  Using (2) we can write $x_i = z_i +
  \sum_{j=1}^{d+1} \alpha_{ij} x'_j$ for each $i$ where $z_i \in
  K_{2}$ and $(\alpha_{ij})$ is an invertible matrix.  Replacing
  $x'_i$ with $x_i$ and $y'_i$ with $\sum_j \alpha_{ij} y'_j$ we may
  assume that $x'_i = x_i$ for each $i$.
  
  Now if $i \leq d$ we have $x_i+y_i \in f_1(1\cc 1) \subset f_2(1\cc
  1)$ by (1) so we may write $x_i+y_i = z'_i + \sum_{j=1}^{d+1}
  \beta_{ij} (x_j+y'_j) = \sum \beta_{ij} x_j + (z'_i + \sum
  \beta_{ij} y'_j)$ where $z'_i \in K_{2}$ and $\beta_{ij} \in \C$.
  Since the last term of this belongs to $K_{2} \oplus \left< y_1,
    \dots, y_{d+1} \right>$, it follows from (4) that $\beta_{ij}$ is
  equal to one if $i=j$ and zero otherwise, so we conclude that
  $x_i+y_i = z'_i + (x_i+y'_i)$.  Thus we have $y_i = y'_i + z'_i$ so
  we may replace $y'_i$ with $y_i$ for $1 \leq i \leq d$.
  
  Finally since $y'_{d+1} \in f_2(0\cc 1)$ we can write $y'_{d+1} =
  z'' + \sum_{j=1}^{d+1} \lambda_j y_j$ by (3) where $z'' \in
  K_{2}$, $\lambda_j \in \C$.  Replacing $y'_{d+1}$ with $y'_{d+1} -
  z''$ we may assume that $z'' = 0$.  Now (1) and (5) imply that
  $f_3(1\cc 1) = K_{3} \oplus \left< x_1+y_1, \dots, x_d+y_d
  \right>$.  Since $x_{d+1}+y'_{d+1} \in \linebreak f_2(1\cc 1)
  \subset f_3(1\cc 1)$ we conclude by (5) and (6) that $\lambda_j = 0$
  for $j \leq d$.  So we have $y'_{d+1} = \lambda_{d+1}\, y_{d+1}$ as
  required.
\end{proof}

\noprt{Spelled wrooooong.}


\section{Proof of the quantum Pieri formula}
\label{sec:proof}

In this section we finally prove Ciocan-Fontanine's quantum Pieri
formula \cite{ciocan-fontanine:on}.  For convenience we will prove the
following equivalent statement of \refthm{thm:partqpieri}.

\newtheorem*{thm1a}{Theorem 1$'$}

\begin{thm1a}
  Let $\alpha = s_r s_{r+1} \cdots s_{a_j}$ and $u,w \in S_n(a)$ be
  permutations, and let $d = (d_1,\dots,d_k)$ be a multidegree, such
  that $\ell(u)+\ell(w)+\ell(\alpha) = \dim \Fl(a;E) + \sum
  (a_{i+1}-a_{i-1})d_i$.  
  The Gromov-Witten invariant $\gw{\pomg_u, \pomg_w, \pomg_\alpha}{d}$
  on $\Fl(a;E)$ is non-zero only if $d$ is a Pieri sequence with
  maximum at position $j$.  In this case we have
\[ \gw{\pomg_u, \pomg_w, \pomg_\alpha}{d} =
   \int_{\Fl(b;E)} \Omega^{(b)}_{\wbar u} 
   \cdot \Omega^{(b)}_{\wbar w}
   \cdot \Omega^{(b)}_{\wbar \alpha}
\]
where $b = a - d = (a_1-d_1,\dots,a_k-d_k)$.
\end{thm1a}

It will be clear from the proof that, if the right hand side of the
identity is non-zero, then $\ell(\wbar u) = \ell(u) - \sum
(a_{i+1}-a_i)d_i = \ell(u)-\ell(\gamma_d)$, so $\wbar u = u \gamma_d$
by \reflemma{lemma:gammabar}, and similarly for $w$.  Therefore the
equivalence with \refthm{thm:partqpieri} is a matter of dualizing the
permutation $w$.  Notice also that the right hand side can only be
equal to zero or one by the classical Pieri formula (\ref{eqn:pieri}).

\begin{proof}
  We first show that if the Gromov-Witten invariant $\gw{\pomg_u,
    \pomg_w, \pomg_\alpha}{d}$ is non-zero, then $d$ is a Pieri
  sequence with maximum at position $j$, and the triple intersection
  on $\Fl(b;E)$ is non-zero as well.  Throughout this proof,
  $F_\bull$, $G_\bull$, and $H_\bull$ will denote full flags of $E$ in
  general position.
  
  Let $C \subset \Fl(a;E)$ be a rational curve of multidegree $d$
  which meets each of the Schubert varieties
  $\Omega_u^{(a)}(F_\bull)$, $\Omega_w^{(a)}(G_\bull)$, and
  $\Omega_\alpha^{(a)}(H_\bull)$.  Let $K_\bull \in \Fl(b;E)$ be the
  kernel of $C$ and set $e_i = a_i - b_i$ for $1 \leq i \leq k$.  Then
  \cite[Lemma~1]{buch:quantum} shows that $e_i \leq d_i$ for all $i$,
  and by \refprop{prop:direct} we have $K_\bull \in \Omega_{\wbar
    u}^{(b)}(F_\bull) \cap \Omega_{\wbar w}^{(b)}(G_\bull) \cap
  \Omega_{\wbar \alpha}^{(b)}(H_\bull)$.  In particular $\ell(\wbar u) +
  \ell(\wbar w) + \ell(\wbar \alpha) \leq \dim \Fl(b;E)$.
By the definition of $\wbar w$ we get $\ell(\wbar \alpha) \geq
\ell(\alpha) - e_j$, $\ell(\wbar u) \geq \ell(u) - \sum_{i=1}^k
(a_{i+1} - a_i) e_i$, and $\ell(\wbar w) \geq \ell(w) - \sum_{i=1}^k
(a_{i+1} - a_i) e_i$.  Thus we obtain
\[ \dim \Fl(b;E) 
   \geq \ell(\wbar u) + \ell(\wbar w) + \ell(\wbar \alpha) 
   \geq \dim \Fl(a;E) - e_j + 
        \sum_{i=1}^k (2 a_i - a_{i-1} - a_{i+1}) e_i \,.
\]
Since $\dim \Fl(b;E) - \dim \Fl(a;E) = \sum (2 a_i - a_{i-1} -
a_{i+1}) e_i + \sum (e_i e_{i+1} - e_i^2)$, this implies that
\[ e_j + \sum e_i e_{i+1} - \sum e_i^2 \geq 0 \,.\]

\reflemma{lemma:pieriseq} therefore shows that $e$ is a Pieri sequence
with maximum at position $j$ and that all the inequalities above must
be satisfied with equality.  In particular we have $d = e$.
Furthermore, since $\ell(\wbar u) + \ell(\wbar w) + \ell(\wbar \alpha)
= \dim \Fl(b;E)$ we must have $\int \Omega_{\wbar u}^{(b)} \cdot
\Omega_{\wbar w}^{(b)} \cdot \Omega_{\wbar \alpha}^{(b)} = 1$ as
required.

On the other hand, if $d$ is a Pieri sequence with maximum at position
$j$, and if the triple intersection on $\Fl(b;E)$ is non-zero, then
the same dimension count shows that $\ell(\wbar u) = \ell(u) -
\ell(\gamma_d)$ and $\ell(\wbar w) = \ell(w) - \ell(\gamma_d)$, so
$\wbar u = u \gamma_d$ and $\wbar w = w \gamma_d$ by
\reflemma{lemma:gammabar}.  Since $\ell(\wbar \alpha) = \ell(\alpha) -
d_j$, we deduce that $\wbar \alpha = s_r s_{r+1} \cdots s_{b_j}$.  We
will continue by showing that the Gromov-Witten invariant
$\gw{\pomg_u, \pomg_w, \pomg_\alpha}{d}$ is non-zero by explicitly
constructing a rational curve $C \subset \Fl(a;E)$ of multidegree $d$,
which meets each of the Schubert varieties $\pomg_u(F_\bull)$,
$\pomg_w(G_\bull)$, and $\pomg_\alpha(H_\bull)$.

As in \S \ref{sec:qpieri} we set $h_p = \min \{i : d_i = p\}$ and $l_p
= \max \{ i : d_i = p \}$ for $1 \leq p \leq d_j$, so that
$\gamma_d(b_{l_p+1}) = a_{h_p}$, while $\gamma_d(i) = i+p$ if $b_{h_p}
< i \leq b_{h_{p+1}}$ or if $b_{l_{p+1}+1} < i < b_{l_p+1}$.  We also
set $E_i = F_{n+1-i} \cap G_i$ for $1 \leq i \leq n$.  Since the flags
are general, it follows that these spaces have dimension one, and $E =
E_1 \oplus \dots \oplus E_n$.

Let $\til w = w_0 {\wbar w} w_b$ be the dual permutation of $\wbar w
\in S_n(b)$.  Then $\wbar u \xra{\wbar \alpha} \til w$.  For each $1
\leq i \leq n$ we define a space $L_i \subset E$ as follows.  If
$\wbar u(i) = \til w(i)$ or if $i > b_j$ we set $L_i = E_{\wbar
  u(i)}$.  Let $B$ be the direct sum of the spaces $L_i$ for which $i
\leq b_j$ and $\wbar u(i) = \til w(i)$.  Notice that $\dim B = r-1$.
When $i \leq b_j$ and $\wbar u(i) \neq \til w(i)$ we then let $L_i$ be
the unique one-dimensional subspace of $E_{\wbar u(i)} \oplus E_{\til
  w(i)}$ such that $(B \oplus L_i) \cap H_{n-r} \neq 0$.  This is well
defined since $B \oplus E_{\wbar u(i)} \oplus E_{\til w(i)}$ has
dimension $r+1$ and since the flags are general, and furthermore we
have $L_i \neq E_{\wbar u(i)}$ and $L_i \neq E_{\til w(i)}$.

It follows from \reflemma{lemma:basis} that the spaces $L_i$ are
linearly independent and that the partial flag $K_\bull \in \Fl(b;E)$
defined by $K_{i} = L_1 \oplus \dots \oplus L_{b_i}$ belongs to the
Schubert variety $\Omega_{\wbar u}^{(b)}(F_\bull)$.  We furthermore
get the same partial flag $K_\bull$ if we take $L_i = E_{\til w(i)}$
for $i > b_j$.  A symmetric argument therefore shows that $K_\bull \in
\Omega_{\wbar w}^{(b)}(G_\bull)$.

Finally, since $(B \oplus L_i) \cap H_{n-r} \neq 0$ for $\ell(\wbar
\alpha)$ different indices $i \leq b_j$, we obtain $\dim(K_{j} \cap
H_{n-r}) = \ell(\wbar \alpha)$ which means that $K_\bull \in \Omega_{\wbar
  \alpha}^{(b)}(H_\bull)$.  By the classical Pieri formula we
therefore conclude that
\begin{equation}
  \Omega_{\wbar u}^{(b)}(F_\bull) \cap \Omega_{\wbar w}^{(b)}(G_\bull)
  \cap \Omega_{\wbar \alpha}^{(b)}(H_\bull) = 
  \{ K_\bull \} \,.
\end{equation}
Notice that this implies that $K_\bull$ must be the kernel of any
rational curve of multidegree $d$ i $\Fl(a;E)$ which passes through
$\pomg_u(F_\bull)$, $\pomg_w(G_\bull)$, and $\pomg_\alpha(H_\bull)$.

Now for each $1 \leq p \leq m = d_j$ choose $x_p \in E_{u(a_{h_p})}$
and $y_p \in E_{w_0 w(a_{h_p})}$ such that $(B \oplus \left<x_p +
  y_p\right>) \cap H_{n-r} \neq 0$.  Notice that since
$\gamma_d(b_{l_p+1}) = a_{h_p}$ we have $u(a_{h_p}) = \wbar
u(b_{l_p+1})$ and $w_0 w(a_{h_p}) = \til w(b_{l_p} + 1)$.  By
\reflemma{lemma:basis} this implies that $B \cap (E_{u(a_{h_p})}
\oplus E_{w_0 w(a_{h_p})}) = 0$, so $x_p$ and $y_p$ can be found.
Since the integers $b_{l_i+1}$ and $b_{l_i}+1$ are all different for
$1 \leq i \leq m$, the same lemma furthermore implies that $x_1,
\dots, x_m, y_1, \dots, y_m$ are linearly independent and that
$K_{{l_p}} \cap \left< x_1, \dots, x_p, y_1, \dots, y_p \right> = 0$
for all $p$.

Let $f : \P^1 \to \Fl(a;E)$ be the morphism which maps a point $(s\cc
t) \in \P^1$ to the partial flag $V_\bull \in \Fl(a;E)$ given by
$V_{i} = K_{i} \oplus \left< s x_1 + t y_1, \dots, s x_{d_i} + t
  y_{d_i} \right>$.  Since $x_p, y_p \in K_{{l_p+1}}$ it follows
that $V_{i} \subset V_{{i+1}}$ for all $i$, so $f$ is well
defined.  Its image $C = f(\P^1) \subset \Fl(a;E)$ is a rational curve
of multidegree $d$.

Notice that $F_{n+1-\wbar u(i)} = F_{n-\wbar u(i)} \oplus L_i$.  If we
set $L'_p = L_{\gamma_d^{-1}(p)}$ for $1 \leq p \leq n$ then
$F_{n+1-u(p)} = F_{n-u(p)} \oplus L'_p$ and the space of dimension
$a_i$ in the partial flag $f(1\cc 0)$ is equal to $L'_1 \oplus \dots
\oplus L'_{a_i}$.  This shows that $f(1\cc 0) \in \pomg_u(F_\bull)$.
A symmetric argument shows that $f(0\cc 1) \in \pomg_w(G_\bull)$.
Finally, since $\dim(K_{j} \cap H_{n-r}) = \ell(\alpha)-d_j$ and since
$(B \oplus \left<x_i+y_i\right>) \cap H_{n-r} \neq 0$ for $1 \leq i
\leq d_j$ we conclude that $f(1\cc 1) \in \pomg_\alpha(H_\bull)$.
This proves that $\gw{\pomg_u, \pomg_w, \pomg_\alpha}{d} \neq 0$.

It remains to be shown that $C$ is the only curve which contributes to
this Gromov-Witten invariant.  Let $f' : \P^1 \to \Fl(a;E)$ be any
rational curve of multidegree $d$ such that $f'(1\cc 0) \in
\pomg_u(F_\bull)$, $f'(0\cc 1) \in \pomg_w(G_\bull)$, and $f'(1\cc 1)
\in \pomg_\alpha(H_\bull)$.  Then the kernel of $f'(\P^1)$ must be
$K_\bull$.  We will show that $f'$ is identical to the map $f$
constructed above.

Set $f'_i = \rho_i \circ f'$ and $f_i = \rho_i \circ f : \P^1 \to
\Gr(a_i,E)$.  We will prove that $f'_i = f_i$ by induction on $i$, the
case $i=0$ being clear.  Assume that $i > 0$ and that $f'_{i-1} =
f_{i-1}$.  If $i \not \in \{ h_1, \dots, h_m \}$ then this follows
because $f'_i(s\cc t) \supset f'_{i-1}(s\cc t) + K_{i} = f_{i-1}(s\cc
t) + K_{i} = f_i(s\cc t)$ by the definition of $f$.

So suppose $i = h_p$ for some $p$.  Then we know that $f'_{h_p-1} =
f_{h_p-1}$.  Since $f'_{l_p+1}(s\cc t)$ must contain
$f_{h_p-1}(s\cc t) + K_{{l_p+1}} = f_{l_p+1}(s\cc t)$ we
furthermore deduce that $f'_{l_p+1} = f_{l_p+1}$.  In particular we
see that the span of $f'_{l_p+1}(\P^1)$ is the space $W =
K_{{l_p+1}} \oplus \left< x_1, \dots, x_{p-1}, y_1, \dots, y_{p-1}
\right>$.

Notice that since $\ell(u \gamma_d) = \ell(u) - \ell(\gamma_d)$ we
have $u(a_{h_\tau}) > u(a_{h_p})$ for $\tau < p$ (see the remarks
after \refthm{thm:partqpieri}.)  Similarly we have $w_0w(a_{h_\tau}) <
w_0w(a_{h_p}) = \til w(b_{l_p}+1) \leq \wbar u(b_{l_p}+1) \leq \wbar
u(b_{l_p+1}) = u(a_{h_p})$.


We claim that the intersection of $W$ with $F_{n+1-u(a_i)}$ is
contained in $K_{i} \oplus \left<x_1, \dots, x_{p-1}\right> \oplus
E_{u(a_i)}$.  To see this, notice at first that $x_1, \dots, x_{p-1}
\in F_{n-u(a_i)}$ and $y_1, \dots, y_{p-1} \in G_{u(a_i)}$ by the
above inequalities.  We will show that for any $\tau \leq b_{l_p+1}$
we have $L_\tau \subset F_{n-u(a_i)}$ or $L_\tau \subset G_{u(a_i)}$.
Again using that $\ell(u \gamma_d) = \ell(u) - \ell(\gamma_d)$ it
follows that $u(a_i) = \wbar u(b_{l_p+1}) \geq \wbar u(\tau)$ for all
$b_{i} < \tau \leq b_{l_p+1}$.  If $b_j < \tau \leq b_{l_p+1}$ we
therefore get $L_\tau = E_{\wbar u(\tau)} \subset G_{u(a_i)}$.  If
$\tau \leq b_j$ then $\wbar u(\tau) \leq \til w(\tau)$.  So if $\wbar
u(\tau) > u(a_i)$ then $L_\tau \subset F_{n-u(a_i)}$.  Finally, if
$\tau \leq b_j$ and $\wbar u(\tau) \leq u(a_i)$ then it follows from
\reflemma{lemma:misc1} that $\til w(\tau) \leq u(a_i)$, so $L_\tau
\subset G_{u(a_i)}$.  Since in particular $L_\tau \subset G_{u(a_i)}$
for all $b_i < \tau \leq b_{l_p+1}$, it follows that $W \cap
F_{n-u(a_i)} \subset K_{i} \oplus \left< x_1,\dots,x_{p-1} \right>$.
Our claim follows from this since $E_{u(a_i)} \subset W$.

Since $F_{n+1-u(\tau)} = F_{n-u(\tau)} \oplus L'_\tau$ for all $\tau$,
it follows that the intersection of $K_{i} \oplus \left<
  x_1,\dots,x_{p-1} \right> = L'_1 \oplus \dots \oplus L'_{a_i-1}$
with $F_{n+1-u(a_i)}$ has dimension $\# \{ \tau \leq a_i \mid u(\tau)
\geq u(a_i) \} - 1$.  Therefore we can write $f'_i(1\cc 0) = K_{i}
\oplus \left< x_1,\dots,x_{p-1},x'_p \right>$ where $x'_p \in
F_{n+1-u(a_i)}$.  Since $x'_p$ must also be contained in $W$, it
follows from the claim that we can take $x'_p \in E_{u(a_i)}$.  A
symmetric argument shows that $f'_i(0\cc 1) = K_{i} \oplus \left< y_1,
  \dots, y_{p-1}, y'_p \right>$ where $y'_p \in E_{w_0 w(a_i)}$.
Applying \reflemma{lemma:classify} to the map $(f_{h_p-1}, f_{h_p},
f_{l_p+1}) : \P^1 \to \Fl(a_{h_p-1}, a_{h_p}, a_{l_p+1}; E)$ we then
conclude that $f'_i(s\cc t) = K_{i} \oplus \left< s x_1 + t y_1,
  \dots, s x_{p-1} + t y_{p-1}, s x'_p + t \lambda y'_p \right>$ for
some $\lambda \in \C^*$.

Finally notice that we can write $f'_j(1\cc 1) = B \oplus \left< x'_p +
  \lambda y'_p \right> \oplus M$ for some subspace $M \subset E$ of
dimension $\ell(\alpha)-1$.  Since the dimension of $f'_j(1\cc 1) \cap
H_{n-r}$ is at least $\ell(\alpha)$, it follows that $(B \oplus \left<
  x'_p + \lambda y'_p \right>) \cap H_{n-r} \neq 0$, so
$\C(x'_p+\lambda y'_p) = \C(x_p+y_p)$.  This shows that $f'_i = f_i$
and finishes the proof.
\end{proof}

\noprt{Spelled wrooooong.}


\providecommand{\bysame}{\leavevmode\hbox to3em{\hrulefill}\thinspace}
\providecommand{\MR}{\relax\ifhmode\unskip\space\fi MR }
\providecommand{\MRhref}[2]{%
  \href{http://www.ams.org/mathscinet-getitem?mr=#1}{#2}
}
\providecommand{\href}[2]{#2}


\end{document}